\documentclass[12pt,twoside]{article}
\usepackage{amsfonts}
\usepackage{amsmath}
\usepackage{amssymb}
\usepackage{latexsym}
\newtheorem{thm}{\bf Theorem}[section]
\newtheorem{cor}[thm]{\bf Corollary}
\newtheorem{lem}[thm]{\bf Lemma}

\newtheorem{exmp}[thm]{\bf Example}

\def\proof{{\parindent0pt {\bf Proof.\ }}}
\def\wdim{{\rm wdim}}
\def\gldim{{\rm glim}}

\def\pd{{\rm pd}}
\def\fd{{\rm fd}}

\def\sup{{\rm sup}}

\newcommand{\cqfd}
{\hspace{1cm}
\rule{2mm}{2mm}%
\medbreak%
\par%
}
\begin{document}
\thispagestyle{empty}


\title{Some homological properties of an amalgamated duplication of a ring along
an ideal}

\author{Mohamed Chhiti and Najib Mahdou\\ \\
Department of Mathematics\\
Faculty of Science and Technology of Fez, P.O. Box 2202\\
University S. M. Ben Abdellah, Fez, Morocco\\
chhiti.med@hotmail.com \\
mahdou@hotmail.com}

\date{}
\maketitle 

\noindent{\large\bf Abstract.}
In this work, we investigate the transfer of some homological properties
from a ring $R$ to his amalgamated duplication along some ideal $I$
of $R$, and then generate new
and original families of rings with these properties.
\bigskip


\begin{section}{Introduction}

Let $R$ be a commutative ring with unit element 1 and let $I$ be a
proper ideal of $R$. The amalgamated duplication of a ring $R$
along an ideal $I$ is a ring that is defined as the following
subring with unit element $(1,1)$ of $R\times R$:
\begin{eqnarray*}
   &R\bowtie I =\{(r,r+i) \mid r\in R,i\in I\}&
\end{eqnarray*}
This construction has been studied, in the general case, and
from the different point of view of pullbacks, by M D'Anna and M
Fontana  \cite{AF2}. Also, M. D'Anna and M. Fontana, in
\cite{AF1}, have considered the case of the amalgamated duplication
of a ring ,in not necessarily Noetherian setting, along a
multiplicative -canonical ideal in the sense
Heinzer-Huckaba-Papick \cite{HHP}. In \cite{A}, M. D'Anna has studied some properties of
$R\bowtie I $, in order to construct reduced Gorenstein rings
assosiated to Cohen-Macaulay rings and has applied this
construction to curve singularities.  On the other
hand, H.R Maimani and S Yassemi, in \cite{MY}, have studied the diameter and
girth of the zero- divisor of the ring $R\bowtie I $. For instance, see \cite{A,AF1,AF2,MY}.

Let $M$ be an $R$-module, the idealization $R\propto M$ (also
called the trivial extention), introduced by Nagata in 1956 (cf
 \cite{N}) is defined as the $R$-module $R\oplus M$ with multiplication
defined by $(r,m)(s,n)=(rs,rn+sm)$. For instance, see
\cite{G,H,KM1,KM2}.

When $I^{2}=0$, the new construction $R\bowtie I $ coincides with
the idealization $R\propto I$. One main difference of this
construction, with respect to idealization is that the ring
$R\bowtie I $ can be a reduced ring (and, in fact, it is always
reduced if $R$ is a domain).

For two rings $A\subset B$, we say that $A$ is a module retract
(or a subring retract) of $B$ if there exists an $A$-module
homomorphism $\varphi : B \rightarrow A$ such that
$\varphi \mid _{A}=id \mid _{A}$. $\varphi$ is called  a
 module retraction map. If such a map $\varphi$ exists, $B$
 contains $A$ as an $A$- module direct summand. We can easily show
 that $R$ is a module retract of $R\bowtie I$, where the module
 retraction map $\varphi$ is defined by $\varphi(r,r+i)=r$.

In this paper, we study the transfer of some homological properties
from a ring $R$ to a ring  $R\bowtie I $. Specially, in section 2,
we prove that $R\bowtie I $ is a Von Neumann
regular ring (resp., a perfect ring) if and only if so is $R$. Also,
we prove that $gldim(R\bowtie I) = \infty $ if $R$ is a domain and $I$
is a principal ideal of $R$. In section 3, we study the coherence of
$R\bowtie I $. More precisely, we prove that if $R$ is a coherent
ring and $I$ is a finitely generated ideal of $R$, then  $R\bowtie
I $ is coherent. And if $I$ contains a regular element, we prove
the converse.

 Recall that if $R$ is a ring and $M$ is an $R$-module, as usual we use
$\pd_R(M)$ and $\fd_R(M)$ to denote the usual projective and flat
dimensions of $M$, respectively. The classical global and weak
dimension of $R$ are respectively denoted by $\gldim(R)$ and $\wdim(R)$. Also,
the Krull dimension of $R$ is denoted by $dim(R)$.
\end{section}

\begin{section}{ Transfer of some homological properties\\ from a ring $R$ to a ring $R\bowtie I$  }

Let $R$ be a commutative ring with identity element 1 and let $I$
be an ideal of $R$. We define $R\bowtie I =\{(r,s)/r,s\in R,s-r\in
I\}$. It is easy to check that $R\bowtie I$ is a subring with unit
element $(1,1)$,  of $R\times R$ (with the usual componentwise
operations) and that $R\bowtie I =\{(r,r+i)/r\in R,i\in I\}$.

It is easy to see that, if $ \pi_{i}(i=1,2)$ are the projections
of $R\times R$ on $R$, then $\pi_{i}(R\bowtie I)=R$ and hence if
$O_{i}=ker(\pi_{i}\backslash R\bowtie I)$. Then $R\bowtie
I/O_{i}\cong R$. Moreover $O_{1}=\{(0,i),i\in I\}$,
 $O_{2}=\{(i,0), i\in I\}$ and $O_{1}\cap O_{2}=(0)$.

 We begin by studying the transfer of Von Neumann regular property.

\begin{thm}\label{1}Let $R$ be a commutative ring  and let  $I$ be a proper
ideal of $R$. Then $R$ is a Von Neumann regular ring if and only
if $R\bowtie I $ is a Von Neumann  regular ring.
\end{thm}

The proof will use the following Lemma.

\begin{lem}\cite[Theorem 3.5]{AF2}\label{4}
\begin{enumerate}
    \item Let $R$  be a commutative ring  and let  $I$ be an ideal
 of $R$. Let $P$ be a prime ideal of $R$ and set:
\begin{eqnarray*}
   &P_{0}=\{(p,p+i)/p\in P,i\in I\cap P\}&  \\
   &P_{1}=\{(p,p+i)/p\in P,i\in I\}&  \\ and
   &P_{2}=\{(p+i,p)/p\in P,i\in I\}&
\end{eqnarray*}
\begin{itemize}
    \item If $I\subseteq P$, then $P_{0}=P_{1}=P_{2}$ is a prime ideal of $R\bowtie I $
    and it is the unique prime ideal of $R\bowtie I $ lying over
    $P$.
    \item If $I\nsubseteq P$, then $P_{1}\neq P_{2}$, $P_{1}\cap
    P_{2}=P_{0}$ and  $P_{1}$  and $P_{2}$ are the
    only prime ideals of $R\bowtie I $ lying over $P$.
\end{itemize}

    \item Let $Q$ be a prime ideal of $R\bowtie I $  and let $O_{1}=\{(0,i) /i\in
    I\}$. Two cases are possible: either $Q \nsupseteq O_{1}$ or  $Q \supseteq
    O_{1}$.

\begin{description}
    \item[a-]If $Q \nsupseteq O_{1}$ , then there exists a unique
    prime ideal $P$ of $R$ $(I\nsubseteq P)$such that
\begin{eqnarray*}
   &Q=P_{2}=\{(p+i,p)/p\in P,i\in I\}&
\end{eqnarray*}
    \item[b-]If $Q \supseteq  O_{1}$, then there exists a unique
    prime ideal $P$ of $R$ such that
\begin{eqnarray*}
   &Q=P_{1}=\{(p,p+i)/p\in P,i\in I\}&
\end{eqnarray*}
\end{description}
\end{enumerate}
\end{lem}

\textbf{Proof of Theorem 2.1}.
Assume that $R$ is a Von Neumann regular ring. Then $R$ is reduced and so
$R\bowtie I$ is reduced by \cite[Theorem 3.5 (a)(vi)]{AF2}.
It remains to show that $dim(R\bowtie I) =0$ by \cite[Remark, p. 5]{H}.

Let $Q$ be a prime ideal of $R\bowtie I$. If $P=Q\cap R$,
then necessarily  $Q\in \{ P_{1},P_{2}\}$ (by Lemma \ref{4}(2)). But
$P$ is a maximal ideal of $R$ since $R$ is a Von Neumann regular
ring. Then $P_{1}$ and $P_{2}$ are maximal ideals of $R\bowtie I$
(by \cite[Theorem 3.5 (a)(vi)]{AF2}). Hence, $Q$ is a maximal
ideal of $R\bowtie I $, as desired.

Conversely, assume that $R\bowtie I $ is a Von Neumann regular ring. By
\cite[Theorem 3.5 (a)(vi)]{AF2}, $R$ is reduced. Let $P$ be a prime
ideal of $R$. By Lemma \ref{4}(1), $P\bowtie I =\{(p,p+i)/p\in P,
i\in I\}$ is a prime ideal of $R\bowtie I$. From \cite[page 7 ]{H}
we get $P\bowtie I$ is a maximal ideal of $R\bowtie I $ and hence
$P$ is a maximal ideal of $R$. Therefore, $dim(R) =0$ and so $R$ is a Von Neumann
regular ring .\cqfd

\begin{cor}Let $R$ be a commutative ring  and let  $I $ be a proper ideal
 of $R$. Then $R$ is a semisimple  ring  if and only if  $R\bowtie I $ is
 a semisimple ring.
\end{cor}

\proof Assume that $R$ be a semisimple ring. Then $R$ is a Noetherian Von Neumann
regular ring. By Theorem \ref{1}, $R\bowtie I $ is
a Von Neumann regular ring  and  by \cite[Corollary 3-3 ]{AF2},
$R\bowtie I $ is Noetherian. Therefore $R\bowtie I $ is
semisimple.

Conversely, assume that $R\bowtie I $ is semisimple. Then  $R\bowtie I $ is
a Noetherian Von Neumann regular ring and so $R$ is a Von
Neumann regular ring (by Theorem \ref{1}) and Noetherian (by
\cite[Corollary 3-3 ]{AF2}). Hence, $R$ is semisimple.\cqfd
\bigskip

A ring $R$ is called a stably coherent ring if for every positive
integer $n$, the polynomial ring in $n$ variables over $R$ is a
coherent ring. Recall that a ring $R$ is  is called a coherent
ring if every finitely generated ideal of $R$ is finitely
presented.

 \begin{cor}Let $R$ be a commutative ring  and let $I $ be a proper ideal
 of $R$. If $R$ is a  Von Neumann regular ring, then $R\bowtie I $
 is a stably coherent ring.
\end{cor}

\proof By Theorem \ref{1} and \cite[Theorem 7.3.1]{G}\cqfd
\bigskip

Now, we are able to construct a new class of non-Noetherian Von
Neumann regular rings.

\begin{exmp}Let $R$ be a non-Noetherian Von Neumann regular ring,
and let $I $ be a proper ideal of $R$. Then, $R\bowtie I $ is
a non-Noetherian Von Neumann regular ring, by
\cite[Corollary 3-3 ]{AF2} and Theorem \ref{1}.
\end{exmp}

We recall that a ring $R$ is called a perfect ring if every flat
$R$-module is a projective $R$-module (see \cite{B}). Secondly,
we study the transfer of perfect property.

\begin{thm}Let $R$ be a commutative ring  and let $I $ be a proper ideal
 of $R$. Then $R$ is a perfect ring if and only if  $R\bowtie I $ is a perfect ring.
\end{thm}

Before proving this Theorem , we need the following Lemmas .

\begin{lem} (\cite[Lemma 2.5.(2)]{M1}\label{2}) \\
Let $(R_{i})_{i=1,2}$ be a family of rings and $E_{i}$
 be an $R_{i}$-module for $i=1,2$. Then $\pd_{R_{1}\times R_{2}}(E_{1}\times
 E_{2})=$\sup$\{\pd_{R_{1}}(E_{1}),\pd_{R_{2}}(E_{2})\}$.

\end{lem}

\begin{lem}\label{3}  Let $(R_{i})_{i=1,2}$ be a family of rings and $E_{i}$
 be an $R_{i}$-module for $i=1,2$. Then $\fd_{R_{1}\times R_{2}}(E_{1}\times
 E_{2})=$\sup$\{\fd_{R_{1}}(E_{1}),\fd_{R_{2}}(E_{2})\}$.
 \end{lem}
 \proof The proof  is analogous to the proof of Lemma \ref{2}.

\begin{lem}\label{5} Let $(R_{i})_{i=1,...,m}$ be a family of rings. Then
$\prod_{i=1}^{m}R_{i}$ is a perfect ring if and only if $R_{i}$ is a
perfect ring for each $i=1,...,m$ .
\end{lem}

\proof The proof is done by induction on $m$ and it suffices to check it
for $m =2$. Let $R_{1}$ and $R_{2}$ be two rings such that
 $R_{1}\times R_{2}$ is perfect. Let $E_{1}$ be a flat $R_{1}$-module
 and let $E_{2}$ be a flat $R_{2}$-module. By Lemma \ref{3}, $E_{1}\times
 E_{2}$ is a flat $R_{1}\times R_{2}$ module and so it is a projective $R_{1}\times
 R_{2}$ module since $ R_{1}\times R_{2}$ is a perfect
 ring. Hence, $E_{1}$ is a projective $R_{1}$-module
 and  $E_{2}$  is a projective  $R_{2}$-module by Lemma \ref{2}; that means that
  $R_{1}$ and $R_{2}$ are perfect rings.

Conversely, assume that $R_{1}$ and $ R_{2}$ are two perfect
rings. Let $E_{1}\times E_{2}$ be a flat $R_{1}\times
R_{2}$-module where $E_{i}$ is an $R_{i} $-module for each $i=1,2
$. By Lemma \ref{3},  $E_{1}$ is a flat $R_{1}$-module and let
$E_{2}$ is a flat $R_{2}$-module; so $E_{1}$ is a projective
$R_{1}$-module and $E_{2}$ is a projective $R_{2}$-module.
Therefore $E_{1}\times E_{2}$ is a projective $R_{1}\times R_{2}$
by Lemma \ref{2}, this means that $R_{1}\times R_{2}$ is a perfect rings.\cqfd

\begin{lem}\label{6}Let $R$ be a commutative ring  and let $I  $ be a proper ideal
 of $R$. Then:
 \begin{enumerate}
    \item An $(R\bowtie I)$- module $M$ is projective if and only if
    $M\otimes_{R\bowtie I}(R\times R)$ is a projective $(R\times R)$
    module and $M/O_{1}M$ is a projective $R$- module.

    \item An $(R\bowtie I)$-module $M$ is flat if and only if
    $M\otimes_{R\bowtie I}(R\times R)$ is a flat $(R\times R)$-module
and $M/O_{1}M$ is a flat $R$- module.
 \end{enumerate}
\end{lem}

\proof Note that  $R\bowtie I$ is a subring of  $R\times R$ and
$O_{1}$ is a common ideal of $R\bowtie I$ and $R\times R$ by
\cite[proposition 3-1]{AF2}. The result follows from \cite[theorem 5-1-1]{G}.\cqfd
\bigskip

\textbf{Proof of Theorem }\ref{2}
Assume that $R$ is a perfect ring and let $M$ be a flat  $(R\bowtie I)$-module.
By Lemma \ref{6}(2), $M\otimes_{R\bowtie I}(R\times R)$ is a
flat $(R\times R)$-module and $M/O_{1}M$ is a flat $R$- module.
Then $M\otimes_{R\bowtie I}R\times R$ is a projective $R\times R$-
module (since $R\times R$ is perfect by Lemma \ref{5}), and $M/O_{1}M$
is a projective $R$-module since $R$ is perfect. By Lemma \ref{6}(1), $M$ is a
projective $(R\bowtie I)$-module and so $R\bowtie I $ is a
perfect ring.

Conversely, assume that $R\bowtie I$ is a perfect ring and let $E$ be a flat $R$- module. Then
$E\otimes_{R}(R\bowtie I)$ is a flat $(R\bowtie I)$-module and so it is a projective $(R\bowtie I)$-module
since $R\bowtie I$ is a perfect ring. In addition, for any $R$-
module $M$ and any $n\geq 1$ we have:
$$Ext^{n}_{R}(E,M\otimes_{R}R\bowtie I)\cong
Ext^{n}_{R}(E\otimes_{R}R\bowtie I,M\otimes_{R}R\bowtie I)$$ (see
\cite[page 118]{CE}) and then $Ext^{n}_{R}(E,M\otimes_{R}R\bowtie
I)=0$. As we note that $M$ is a direct summand of
$M\otimes_{R}R\bowtie I$ since $R$ is a module retract of
$R\bowtie I$, $Ext^{n}_{R}(E,M)=0$ for all  $n\geq 1$ and all $R$-
module $M$. This means that $E$ is a projective $R$- module
and so $R$ is a perfect ring.\cqfd \bigskip

We say that a ring
$R$ is Steinitz if any linearly independent subset of a free
$R$-module $F$ can be extended to a basis of $F$ by adjoining
element of a given basis. In \cite[proposition 5.4]{CP}, Cox and
Pendleton showed that Steinitz rings are precisely the perfect
local rings.

By the above Theorem and since $R\bowtie I$ is local if and only if $R$ is local,
we obtain:

\begin{cor}
Let $R$ be a commutative ring  and let $I$ be a proper ideal
 of $R$. Then $R$ is a Steinitz ring if and only if  $R\bowtie I $ is a Steinitz ring.
\end{cor}

\begin{exmp}Let $R=K[X]/(X^{2})$ where $K$ is a field and $X$ an
indetrminate. Then $(K[X]/(X^{2}))\bowtie(\overline{X})$ is a
Steinitz ring.
\end{exmp}

For a nonnegative integer $n$, an $R$-module $E$ is $n$-presented
if there is an exact sequence $F_n\rightarrow F_{n-1} \rightarrow
...\rightarrow F_0 \rightarrow E\rightarrow 0$ in which each $F_i$
is a finitely generated free $R$-module. In
particular, ``$0$-presented'' means finitely generated and
``$1$-presented'' means finitely presented.  \\
Given nonnegative integers $n$ and $d$, a ring $R$ is called an
$(n,d)$-ring if every $n$-presented $R$-module has projective
dimension $\leq d$ ; and $R$ is called a weak $(n,d)$-ring  if every
$n$-presented cyclic $R$-module has projective dimension $\leq d$
(equivalently, if every $(n-1)$-presented ideal of $R$ has
projective dimension $\leq d-1$). For instance, the
$(0,1)$-domains are the Dedekind domains, the $(1,1)$-domains are
the Pr\"ufer domains, and the $(1,0)$-rings are the von Neumann
regular rings. See for instance (\cite{C}, \cite{KM1}, \cite{KM2}, \cite{M1}, \cite{M2}). \\
Now, we give a wide class of rings which are not a weak $(n,d)$-ring (and so not an $(n,d)$-ring)
 for each positive integers $n$ and $d$.

\begin{thm}\label{8} Let $R$ be an integral domain and let  $I (\not= 0)$ be a principal ideal of
$R$. Then $R$ is not a weak $(n,d)$-ring (and so not an $(n,d)$-ring) for each positive integers $n$ and $d$. In particular, $wdim(R\bowtie I) =gldim(R\bowtie I)=\infty$.
\end{thm}

Before proving this Theorem, we need the following Lemma.

\begin{lem}\label{7} Let $R$ be a commutative ring  and let $I (\not= 0)$ be a principal ideal of
$R$, then  $O_{1}=\{(0,i),i\in I\}$ and $O_{2}=\{(i,0), i\in I\}$
are principal ideals of $R\bowtie I $.
\end{lem}

\proof Let $(0,i)$ be an element of $O_{1}$. Since $I$ is a principal
ideal of $R$, then there exists $a\in I$ such that $I =Ra$ and so
  $(0,i)=(0,ra)=(r+j,r)(0,a)$ for some $r \in R$ and for all
$j\in I$. Hence, $O_{1}$ is a principal ideal of $R\bowtie I $ generated by $(0,a)$.
Also, $O_{2}$ is a principal ideal generated by
$(a,0)$ by the same argument, as desired.\cqfd\bigskip

\textbf{Proof of Theorem }\ref{8}.
Let $a\in I$ such that $I =Ra$. By lemma  \ref{7}, $O_{1}$ and $O_{2}$
are principal  ideals of $R\bowtie I $.  Consider the short
exact sequence of $R\bowtie I $-modules:
$$(1)\qquad0\rightarrow ker (u)\rightarrow R\bowtie I\stackrel{u}\rightarrow
O_{1}\rightarrow 0$$

where $u(r,r+i)=(r,r+i)(0,a) =(0,(r+i)a)$. Then,
$ker(u)=\{(r,0)\in R\bowtie I /r\in I\} =O_{2}$. Consider the
short exact sequence of $R\bowtie I $-modules:
$$(2)\qquad0\rightarrow ker (v)\rightarrow R\bowtie I\stackrel{v}\rightarrow
O_{2}\rightarrow 0$$

where $v(r,r+i)=(r,r+i)(a,0) =(ra,0)$. Then,
$ker(v)=\{(0,i)\in R\bowtie I /i\in I\}=O_{1}$. Therefore,
$O_{1}$ (resp., $O_{2}$) is $m$-presented for each positive integer $m$ by the above two exact sequences. It remains to
show that $pd_{R\bowtie I}(O_{1}) = \infty $ (or $pd_{R\bowtie I}(O_{2}) = \infty $).

We claim that $O_{1}$ and $O_{2}$ are not projective. Deny. Then $O_{1}$ is
projective and so the short exact sequence (1) splits. Then
$O_{2}$ is generated by an idempotent element $(x,0)$, such that
$x (\neq 0) \in I$. Hence, $(x,0)^{2}=(x,0)(x,0)=(x^{2},0)=(x,0)$, then
$x^{2}=x$, and so $x=1$ or $x=0$, a contradiction (since $x\in I$
and $x\neq0$). Therefore, $O_{1}$ is not projective. Similar
arguments show that $O_{2}$ is not projective. A combination of (1) and
(2) yields $pd_{R\bowtie I}(O_{1})=pd_{R\bowtie
I}(O_{2})+1$ and $pd_{R\bowtie I}(O_{2})=pd_{R\bowtie I}(O_{1})+1$
then, $pd_{R\bowtie I}(O_{1})=pd_{R\bowtie
I}(O_{2})+1+1=pd_{R\bowtie I}(O_{1})+2$. Consequently,
the projective dimension of $O_{1}$ (resp., $O_{1}$) has to be infinite, as desired. \cqfd \bigskip

If $R$ is a principal domain, we obtain:

\begin{cor}Let $R$ be a principal domain and let  $I$ be a proper ideal
 of $R$. Then $R$ is not a weak $(n,d)$-ring (and so not an $(n,d)$-ring) for each positive integers $n$ and $d$. In particular, $wdim(R\bowtie I) =gldim(R\bowtie I)=\infty$.
\end{cor}
\end{section}

\begin{section}{The coherence of $R\bowtie I $}

An $R$- module $M$ is called a  coherent $R$ module, if it is a
finitely generated  and every finitely generated submodule of $M$
is finitely presented.

A ring $R$ is called a coherent ring if it is a coherent
module over itself, that is, if every finitely generated ideal of
$R$ is finitely presented, equivalently, if $(0:a)$ and $I\cap J$
are finitely generated for every $a\in R$ and any two finitely
generated ideals $I$ and  $J$ of $R$ (by \cite[Theorem 2.2.3]{G}).
Examples of coherent rings are Noetherian rings, Boolean algebras,
Von Neumann regular rings, valuation rings, and Pr\"{u}fer/
semihereditary rings. See for instance [\cite{G}].

\begin{thm}\label{9}Let $R$ be a commutative ring  and let  $I$ be a proper ideal
 of $R$. Then:
\begin{enumerate}
    \item If  $R\bowtie I $ is coherent, then $R$ is coherent .
    \item If $R$ is a coherent ring and $I$ is a finitely generated ideal of
$R$, then  $R\bowtie I $ is coherent.
    \item  Assume that $I$ contains a regular element. Then $R\bowtie I $ is
a coherent ring if and only if $R$ is a coherent ring and $I$ is a
finitely generated ideal of $R$.

\end{enumerate}
\end{thm}

We need the following Lemma before proving this Theorem.

\begin{lem} (\cite[Theorem 2.4.1]{G}\label{10}).
Let $R$ be a commutative  ring and
let $I$ be a proper ideal of $R$, then :
\begin{enumerate}
\item If $R$ is a coherent ring and $I$ is a finitely generated ideal of
$R$, then $R/ I$ is a coherent ring.
\item If  $R/ I$ is a coherent ring and $I$ is a coherent
$R$ module, then $R$ is a coherent ring.
\end{enumerate}
\end{lem}

\textbf{Proof of Theorem }\ref{9}.
\begin{enumerate}
    \item Let $L=\sum_{i=1}^{n}Ra_{i}$ be a finitely generated ideal
of $R$, and set $J :=\sum_{i=1}^{n}(R\bowtie I)(a_{i},a_{i})$.
Consider the exact sequence of $R\bowtie I $  -modules:
$$0\rightarrow ker (u)\rightarrow (R\bowtie
I)^{n}\stackrel{u}\rightarrow J\rightarrow 0$$

where
$u(r_{i},r_{i}+e_{i})_{1\leq i\leq
n}=\sum_{i=1}^{n}(r_{i},r_{i}+e_{i})(a_{i},a_{i}) =(\sum_{i=1}^{n}a_{i}r_{i},
   \sum_{i=1}^{n}a_{i}r_{i}+\sum_{i=1}^{n}a_{i}e_{i})$. Thus
$ker(u)=\{(r_{i},r_{i}+e_{i})_{1\leq i\leq n}\in(R\bowtie I)^{n}/
\sum_{i=1}^{n}r_{i}a_{i}=0 ,\sum_{i=1}^{n}e_{i}a_{i} =0\}$. On
other hand , consider the exact sequence of $R$-modules:
$$0\rightarrow ker( v)\rightarrow R^{n}\stackrel{v}\rightarrow
L\rightarrow 0$$

where $v(b_{i})=\sum_{i=1}^{n}b_{i}a_{i}$. Hence,
$ker( u)=\{(r_{i},r_{i}+e_{i})_{1\leq i\leq n}\in(R\bowtie
I)^{n}/r_{i}\in ker( v);e_{i}\in I^{n}\cap ker( v)\} $. But $J$ is
a finitely presented since it is finitely generated and $R\bowtie I $ is coherent. Hence, $ker (u)
$ is a finitely generated $(R\bowtie I)$-module and so $ker(
v)$ is a finitely generated $R$ -module . Therefore, $L$ is a
finitely presented ideal of $R$ and so $R$ is coherent.

   \item Since $I$ is a finitely generated ideal of $R$,
then $O_{1}$ and $O_{2}$ are a finitely generated ideals of
$R\bowtie I $. Hence, $R\bowtie I $ is a coherent ring by Lemma
\ref{10} since $R$ is a coherent ring and $R\bowtie I
/O_{i}\cong R$, as desired.

    \item Assume that $R\bowtie I $ is a coherent ring. Then $R$ is a
coherent ring by 1). Now, we prove that $I$ is a finitely generated
ideal of $R$. Let $m$ be a non zero element of $I$ and set
$c=(m,0)\in R\bowtie I$. Then:

\begin{eqnarray*}
   (0:c) & =  & \{(r,r+i)\in  R\bowtie I /(r,r+i)(m,0)=0\} \\
    & =   &\{(r,r+i)\in  R\bowtie I /rm=0\}  \\
   & =  &\{(r,r+i)\in  R\bowtie I /r=0\}  \\
   & =  &\{(0,i)\in  R\bowtie I /i\in I\}  \\
   & =  &O_{1}.
\end{eqnarray*}

Since $R\bowtie I $ is a coherent ring, then $(0:c)$ is a finitely
generated ideal of $R\bowtie I $ and so $O_{1}$ is a finitely
generated ideal of  $R\bowtie I $. This means that $I$ is a finitely
generated ideal of $R$.
\\Conversely if $R$ is a coherent ring and  $I$ is a finitely generated ideal of
$R$, then $R\bowtie I $ is a coherent ring by Lemma
\ref{10}(2) and this completes the proof of Theorem 3.1.\cqfd\bigskip
\end{enumerate}

If $R$ is an integral domain, we obtain:

\begin{cor} Let $R$ be an integral domain and let  $I$ be a proper ideal
 of $R$ . Then $R\bowtie I $ is
a coherent ring if and only if $R$ is a coherent ring and $I$ is a
finitely generated ideal of $R$.
\end{cor}

In general, $R\bowtie I $ is a coherent ring doesn't imply that $I$
is a finitely generated of $R$ as the following example shows:

\begin{exmp} Let $R$ be a  non-Noetherian Von Neumann regular ring and let
$I$ be a non finitely generated ideal of $R$ (see for example \cite {C}).
Then $R\bowtie I $ is a coherent ring but $I$ is not a
finitely generated.
\end{exmp}

Now, we are able to construct a new class of non-Noetherian rings.

\begin{exmp} Let $R$ be a non-Noetherian coherent ring and let $I$ be a finitely generated ideal of $R$. Then:
\begin{enumerate}
\item $R\bowtie I $  is a coherent ring by Theorem 3.1(2).
\item $R\bowtie I $ is non-Noetherian by \cite[Corollary 3.3]{AF2} since $R$ is non-Noetherian.
\end{enumerate}
\end{exmp}

We recall that an $R$- module $M$ is called a uniformly coherent
$R$ module, if $M$ is a finitely generated $R$ module and there is
a map $\phi :\mathbb{N}\rightarrow\mathbb{N}$, where $\mathbb{N}$
denotes the natural numbers, such that for every $n\in
\mathbb{N}$, and any nonzero homomorphism $f: R^{n}\rightarrow M$,
$ker(f)$ can be generated by $\phi(n)$ elements.

A ring $R$ is called an uniformly coherent ring if $R$ is uniformly
coherent as a module over itself.

Recall that an uniformly coherent is a coherent ring (by \cite[Theorem 6.1.1]{G}).
Also, there exists Noetherian rings which are not uniformly coherent (see \cite[p. 191]{G}).
See for instance [ \cite[Chapter 6]{G}].

\begin{thm}\label{11} Let $R$ be a Noetherian ring  and let  $I$ be a nilpotent ideal
 of $R$. Then $R$ is an  uniformly coherent ring if and only if $R\bowtie I $
 is an uniformly coherent ring.
\end{thm}

We need the following Lemma before proving this Theorem.

\begin{lem}\label{12} Let $R$ be a commutative  ring and let $I$ be  a finitely
generated ideal of $R$. If $R\bowtie I $ is an  uniformly
coherent ring then so is $R$.
\end{lem}

\proof  The ideal $O_{1} :=\{(0,i),i\in I\}$ is a finitely generated ideal
of $R\bowtie I$ since $I$ is a finitely generated ideal of
$R$. Hence,  $R :=\cong R\bowtie I
/O_{1}$ is an uniformly coherent ring by \cite[Corollary6-1-6]{G}, as desired. \cqfd\bigskip

\textbf{Proof of Theorem }\ref{11}.

If $R\bowtie I$ is an uniformly coherent ring, then so is $R$ by Lemma \ref{12}
since $R$ is Noetherian. Conversely, assume that $R$ is an
uniformly coherent ring.
Let $\varphi: R\bowtie I\rightarrow (R\bowtie I) /O_{1}$
 be a ring epimorphism. Since $R\bowtie I$ is Noetherian (since $R$ is Noetherian), then
 $(R\bowtie I) /O_{1}$ is a finitely presented  $R\bowtie I$
 module. On other hand, $O_{1}$ is nilpotent (since $I$ is
 nilpotent), and $(R\bowtie I) /O_{1} (:\cong R)$ is  uniformly coherent. Hence,
$R\bowtie I $ is an  uniformly coherent ring by \cite[theorem 6-1-8]{G}.
\end{section}





\begin{thebibliography}{999}\addcontentsline{toc}{section}{\protect\numberline{}{Bibliography}}

      \bibitem{A} M. D'Anna , A construction of Gorenstein rings, J. Algebra {\bf 306} (2006), no. 2, 507-519.

      \bibitem{AF1} M. D'Anna and M. Fontana, The amalgamated duplication
       of a ring along a multiplicative-canonical ideal, Ark. Mat. {\bf 45} (2007), no. 2, 241-252.

      \bibitem{AF2} M. D'Anna and M. Fontana, An amalgamated duplication
      of a ring along an ideal: the basic properties. J. Algebra Appl. {\bf 6} (2007), no. 3, 443-459.

      \bibitem{B} H. Bass, Finitistic dimension and a homological
      generalization of semi-primary rings, Trans. Amer. Math.
      Soc. 95 (1960), 466-488.

      \bibitem{C} D. L. Costa, Parameterising families of
      non-Noetherian rings, Comm. Algebra, 22(1994), 3997-4011

      \bibitem{CE} H. Cartan and S. Eilenberg, Homological Algebra,
      Princeton Univ. Press. Princeton (1956).

      \bibitem{CP} S. H. Cox and R.L. Pendleton, Rings for which
      certain flat modules are projective, Amer.Maths. Soc, 150(1970),
      139-156.
      \bibitem{G} S. Glaz; Commutative Coherent Rings,
     Springer-Verlag, Lecture Notes in Mathematics, 1371 (1989).

      \bibitem{H} J. A. Huckaba, Commutative Coherent Rings with Zero
     Divizors. Marcel Dekker, New York Basel, (1988).

      \bibitem{HHP}W. Heinzer, J. Huckaba and I. Papick, m-canonical
      ideals in integral domains, Comm.Algebra 26(1998), 3021-3043.


     \bibitem{KM1} S. Kabbaj and N. Mahdou ; Trivial Extensions Defined by coherent-like condition,
      comm.Algebra . Marcel Dekker (2004) ,3937-3953

      \bibitem{KM2} S. Kabbaj and N. Mahdou ; Trivial extensions of local rings and a conjecture of Costa,
         Lecture  Notes in Pure and Appl. Math., Vol.231, Marcel Dekker, New York, (2003), 301-312.

     \bibitem{M1} N. Mahdou; On Costa's conjecture, Comm.
     Algebra, 29 (7) (2001), 2775-2785.

     \bibitem{M2} N. Mahdou; On 2-Von Neumann regular rings, Comm. Algebra 33 (10) (2005), 3489-3496.

      \bibitem{MY} H. R. Maimani and S. Yassemi, Zero-divisor graphs of amalgamated duplication
      of a ring along an ideal. J. Pure Appl. Algebra {\bf 212} (2008), no. 1, 168-174.

      \bibitem{N} M. Nagata, Local Rings, Interscience, New york, (1962).

      \bibitem{NN}B. Nashier and W. Nichols, On Steinitz properties, Arch. Math, 57 (1991), 247-253.

      \bibitem{Ro} J. Rotman; An introduction to homological
      algebra, Press, New York, 25(1979).




\end{thebibliography}
\end{document}